\documentclass[a4paper,twoside,11pt]{article}

\usepackage{graphicx,boxedminipage,latexsym}
\usepackage[dvips]{color}

\usepackage[english]{babel}

\usepackage{amsmath,amssymb}

\usepackage{amsthm}

\usepackage{amscd}

\usepackage{amsfonts}

\usepackage{mathrsfs}

\begin{document}

\title{Linear stability of projected canonical curves with applications to the slope of fibred surfaces}
\author{M. A. Barja\footnote{Partially supported by DGICYT BFM2003-06001
(Ministerio de Educaci\'on y Ciencia) and by 2005SGR00557
(Generalitat de Catalunya)}, L. Stoppino\footnote{Partially
supported by PRIN 2005 ``Spazi di Moduli e Teorie di Lie''. }}
\date{}

\newcommand{\Q}{\mathbb{Q}}
\newcommand{\Z}{\mathbb{Z}}
\newcommand{\R}{\mathcal{R}}
\newcommand{\pr}{\mathbb P}
\newcommand{\sym}{\mbox{\upshape{Sym}}}
\newcommand{\of}{\omega_f}
\newcommand{\rk}{\mbox{\upshape{rank}}}
\newcommand{\av}{``}
\newcommand{\om}{\omega}
\newcommand{\og}{\omega_\beta}
\newcommand{\oa}{\omega_\alpha}
\newcommand{\op}{\omega_{f'}}
\newcommand{\g}{\gamma}
\newcommand{\y}{\Gamma}
\newcommand{\el}{\mathcal L}
\newcommand{\oo}{\mathcal O}
\newcommand{\surjarrow}{\rightarrow\!\!\!\!\rightarrow}
\newcommand{\rd}{\mbox{\upshape{red.deg}}}
\newcommand{\kap}{\mathbb{C}}
\newcommand{\pil}{\pi_{\Lambda}}
\newcommand{\cliff}{\mbox{\upshape{Cliff}}}
\newcommand{\pic}{\mbox{\upshape{Pic}}}
\newcommand{\fidi}{\varphi_{|D|}}
\newcommand{\gon}{\mbox{\upshape{gon}}}
\newcommand{\ci}{c_1}
\newcommand{\ch}{\mbox{\upshape ch}}
\newcommand{\td}{\mbox{\upshape td}}
\newcommand{\kapp}{\mbox{\itshape {\textbf k}}}
\newcommand{\Si}{\Sigma}
\newcommand{\sig}{\Sigma}
\newcommand{\ld}{\Lambda}
\newcommand{\qb}{{q_f}}
\newcommand{\spn}{\mbox{\upshape span}}
\newcommand{\ann}{\mbox{\upshape Ann}}
\newcommand{\gr}{\mbox{\upshape Gr}}
\newcommand{\xt}{{X_t}}

\newtheorem{teo}{Theorem}[section]
\newtheorem{prop}[teo]{Proposition}
\newtheorem{conjecture}[teo]{Conjecture}
\newtheorem{lem}[teo]{Lemma}
\newtheorem{cor}[teo]{Corollary}
\newtheorem{rem}[teo]{Remark}
\newtheorem{defi}[teo]{Definition}
\newtheorem{ex}[teo]{Example}
\newtheorem{ass}[teo]{Assumption}

\renewcommand{\theequation}{\arabic{section}.\arabic{equation}}

\pagestyle{myheadings} \markboth{\small{M. A. Barja, L.
Stoppino}}{\small{\textit{Linear stability of projected canonical
curves...}}}

\maketitle

\begin{abstract}
Let $f\colon S\longrightarrow B$ be a non locally trivial  relatively minimal fibred surface.
We prove a lower bound for the slope of $f$ depending increasingly from the
relative irregularity of $f$ and the Clifford index of the general fibres.
\end{abstract}

\bigskip

{\scriptsize {\it Key Words and Prhases.} Fibrarion, Slope,
Relative irregularity, Clifford index.}



\section{Introduction and preliminaries}

Let $f\colon S \longrightarrow B$ be a surjective holomorphic map with
connected fibres from a complex smooth projective surface $S$
onto a complex smooth curve $B$. We always assume that it is
relatively minimal, i.e., that there is no $(-1)-$rational curve
contained in a fibre of $f$.
Let $F$ be a general fibre.
We call $f$ a {\it fibration} of
genus $g$ whenever $g=g(F)$; we also set $b=g(B)$. The fibration
is called {\it smooth} if all its fibres are smooth, {\it
isotrivial} if all its smooth fibres are reciprocally isomorphic,
and {\it locally trivial} if it is smooth and isotrivial (i.e. an holomorphic fibre bundle).

Let $\omega _S$ be the canonical line bundle of $S$ and $K_S$ any
canonical divisor. Set $p_g=h^0(S,\omega_S)$, $q=h^1(S,\omega_S)$,
$\chi {\cal O}_S=p_g-q+1$ and let $e(X)$ be the topological Euler
characteristic of X. We consider the following relative
invariants:
\begin{eqnarray}
&&K^{2}_f=(K_{S}-f^{\ast}K_{B})^{2}=K^{2}_{S}-8(b-1)(g-1)\nonumber\\
&&\chi_{f}=\mbox {deg}f_{\ast}\omega_{S/B}=
\chi {\cal O}_{S}-(b-1)(g-1)\nonumber\\
&&e_{f}=e(S)-e(B)e(F)=e(S)-4(b-1)(g-1)\nonumber\\
&&q_f=q(S)-b.\nonumber
\end{eqnarray}

We have the following classical results, when $g \geq 2$:

\begin{itemize}

\item[(i)] (Noether) $12 \chi_f=e_f+K^2_f.$

\item[(ii)] (Zeuthen-Segre) $e_f \geq 0$. Moreover, $e_f=0$ if and only
if $f$ is smooth.

\item[(iii)] (Arakelov) $K^2_f \geq 0$. Moreover, if $K^2_f=0$
then $f$ is isotrivial.

\item[(iv)] $\chi_f \geq 0$. Moreover, $\chi_f=0$ if and
only if $f$ is locally trivial.

\item[(v)] $0 \leq q_f \leq g $. When $b\geq 1$ $q_f=0$ if and
only if $f$ is the Albanese map of $S$. On the other hand, $q_f=g$
if and only if  $S=B \times F$ (cf. \cite{Bv}).

\end{itemize}

We say that $f$ is a {\it non-Albanese fibration} if $q_f>0$.

\medskip

When $f$ is not locally trivial, Xiao (cf. \cite{X}) defines the
slope of $f$ as

$$s(f)=\frac{K^2_f}{\chi_f}.$$

It follows immediately from Noether's equality that $0 \leq s(f)
\leq 12$.

\medskip

We are mostly concerned with a lower bound of the slope. The main
known result in this direction is:

\bigskip

{\it If $g \geq 2$ and $f$ is not locally trivial, then $s(f) \geq
4-\frac{4}{g}$}.

\bigskip

\noindent which is known as the {\it slope inequality}. It was
first proven by Horikawa and Persson for hyperelliptic fibrations.
Xiao gives a proof for general fibrations (cf. \cite{X}) and,
independently, Cornalba and Harris prove it for {\it semistable}
fibrations (i.e., for fibrations where {\it all} the fibres are
semistable curves in the sense of Deligne and Mumford).
Later on, in \cite{LS} it has been proved a generalization of their
method which can be applied to any fibration.

The slope of a fibration turns out to be sensible to a lot of
geometric properties, both of the fibres of $f$ and of $S$ (see
\cite{A-K} for a complete reference).

We like to pay attention to the influence of the relative
irregularity of the the fibration, $q_f$.
In view of our argument, also
the Clifford index of $f$ appears closely related to this
problem. In \cite{konnocliff}, there is a very interesting attempt
to exhibit the lower bound of the slope as an increasing function
of the Clifford index. This seems very clear
for hyperelliptic or trigonal fibrations (see also\cite{konnotrig}, \cite {S-F}), and
for general Clifford index, but in the intermediate cases,
generality conditions are necessary (\cite{B}).

In the case of the relative irregularity $q_f$, it seems again
that the lower bound of the slope should be an increasing function
of $q_f$.
A crucial point where the relative irregularity $q_f$
appears in a fibration is given by the so called Fujita
decomposition:

$$f_*\omega _f={\cal A} \oplus {\cal O}_B^{\oplus q_f}$$

\noindent which also produces a decomposition of the relative
Jacobian fibration associated to $f$.
In particular, notice that  a general fibre of   a non-Albanese fibration has
 non simple Jacobian.

The first result which manifests the influence of $q_f$ on the
slope is due to Xiao (\cite{X}): $s(f)\geq 4$ whenever $q_f>0$ and
equality holds only if $q_f=1$.
Explicit lower bounds depending on
$q_f$ are given in \cite{KoIrr} and in \cite{B-Z}, but they are rather complicated and
seem far to be sharp.
However, from these results it seems clear that there should be a lower bound for the slope
which is an increasing function of the relative irregularity.

\bigskip

We conjecture the following simple behavior for the bound.

\begin{conjecture}\label{conjecture}
Let $f\colon S \longrightarrow B$ be a fibration of genus $g$, with
relative irregularity  $q_f <g-1$. Then $$s(f) \geq
4\frac{g-1}{g-q_f}.$$
\end{conjecture}
This bound, if true, is sharp (Example \ref{doppie}).
Apart from the aforementioned result of Xiao, some other evidences for this conjecture are
the following.

\begin{itemize}

\item It is true when $\mathbb{P}({\cal O}_B^{\oplus q_f})$ does not
meet the general fibre and the projection from it induces a
birational and linearly stable map (Remark \ref{se...}).
\item It is true when $\mathbb{P}({\cal O}_B^{\oplus q_f})$ does not
meet the general fibre and $\cal A$ is a semistable sheaf on $B$ (Remark \ref{semistable}).
\item There is an analogous canonical decomposition of $f_*\omega _f$ in case the fibration
is a {\it double cover fibration}. In that situation, the corresponding
conjectured bound holds (see Example \ref{doppie} and \cite{CS}).
\item In a semistable fibration with $s$ singular fibres, Vojta
proves the following inequality $$K^2_f \leq (2g-2)(2b-2+s)$$
\noindent which combined with slope inequality gives
$${\chi _f}\leq \frac{g}{2}(2b-2+s).$$
However, a sharper bound of this type holds (cf. \cite{Ara} and
\cite{V-Z}), namely
$${\chi _f} \leq \frac{q-q_f}{2}(2b-2+s)$$
\noindent which is exactly the bound we would obtain using our
conjectured bound instead of slope inequality in Vojta formula.

\end{itemize}

To our knowledge, the only known  counterexamples  to the bound
above belong to the extremal case $q_f=g-1$ (cf. \cite{Pir92} and
Remark \ref{controesempio}).

\bigskip

Our approach is the following. Consider any vector subbundle
${\cal F} \subseteq f_* \omega _f$. The inclusion induces a linear
system on $F$ which is just the projection
$$\mathbb{P}(f_* \omega _f) \longrightarrow \mathbb{P}(\cal F)$$
\noindent restricted to the canonical embedding of $F$ (assume it
is non hyperelliptic). Information about the degree and rank of
this linear system is the main ingredient for applying Xiao's
method. In some cases this information allows to conclude that the
projection is linearly stable; roughly speaking, this means that
any linear subsystem can only increase the ratio between the
degree and the rank (see section 2 for a more precise definition).
In the case of curves linear stability implies Hilbert stability
and so we can also apply Cornalba-Harris method to study a lower
bound of the slope $s(f)$.

With this purpose, we start in section 2 studying when a
projection of a canonical curve is linearly stable. Our main
result in this direction is

\begin{teo}
Let $C\subset \pr^{g-1}$ be a canonical non-hyperelliptic curve. Let $\sig\subset
\pr^{g-1}$ be a $(s-1)$-space. Let $k$ be any
positive integer smaller or equal to $\min\left\{[s/2], [\cliff (C)/2]\right\}$.
Then there is a non-empty open set of $(k-1)$-spaces contained in
$\sig$ that induce linearly stable projections.
\end{teo}

In section 3 we use this information to study a lower bound of the
slope of non-Albanese fibrations. We obtain

\begin{teo}\label{resultado}
Let $f\colon X\rightarrow B$ be a fibred surface. Let $m:=\min
\{\cliff(f),\qb\}$. Then the slope of $f$ satisfies the inequality
$$
s(f)\geq 4\frac{g-1}{g-[m/2]}.
$$
\end{teo}

Although the main ingredient for the theorem is the result of
linear stability of section 2, we give two different proofs of
this result, one by applying  Xiao's method and another
one using the one of Cornalba-Harris. We present this fact as another instance
that, at least in the case of surfaces, both methods, of different
nature, produce the same results. We believe that this parallelism
(which does not clearly hold for higher dimensions) merits further
investigation.

\medskip

The two invariants involved in our main result, the relative irregularity and the Clifford index,
are of very different nature. Theorem \ref{resultado} gives  a strong inequality for big values
of {\em both} these invariants.
It is therefore important to verify that this two quantities are independent,
and in particular that they can grow simultaneously.
In section 4 we provide examples of fibred surfaces with both $\qb$ and
$\cliff (f)$ large, but also of surfaces with large $\qb$ and small $\cliff (f)$, and vice versa.

\medskip

\noindent{\bf Acknowledgments} We thank Maurizio Cornalba,
who gave us uncountably helpful suggestions.
We would also like to thank Gian Pietro Pirola, for many useful discussions, and
Andreas Leopold Knutsen for his kind help with the examples of the last section.

\section{Linear stability of projections of a canonical curve}
In this section we prove, under suitable assumptions, the linear stability of general
projections of canonical curves. This is the key property that allows us to apply both Xiao's and
Cornalba-Harris method in the second part of the paper.


The notion of linear stability  was first defined by Mumford in
\cite{Mum} for embeddings in projective spaces.
The following is a natural generalization for curves with any map to projective spaces.
For a more general treatment, see \cite{LS}.

Let  $C$ be a smooth curve, together with a non-degenerate map in a projective space
$\psi\colon C\rightarrow \mathbb P^{r}$.
Consider  the base point free linear series $\mathcal A$ associated to the morphism $\widetilde\psi$
obtained eliminating the base points of $\psi$.
If $d$ is the degree of $\mathcal A$ and $r$ is its dimension (i.e. $\mathcal A$ is a $g^r_d$),
we define the {\em reduced degree} of the pair $(C,\psi)$ as
$$
\rd(C,\psi)\colon= \frac{d}{r}
$$
(we will also use the notation $\rd (C,\mathcal A)$, or  $\rd(C,V)$, where
$V\subseteq H^0(C,\psi^*\oo_{\mathbb P^r}(1))$ is the linear system such that $\mathcal A=\pr (V)$).

\begin{defi}
With the same notations as above,
we say that $\psi\colon C\rightarrow \pr^r$ is linearly semistable (resp. stable) if for any projection
$\pi$ on a positive dimensional projective space,
 $$\rd (C,\pi\circ\psi)\geq \rd(C,\psi)$$ (resp. $\rd (C,\pi\circ\psi)> \rd(C,\psi)$).
\end{defi}
In other words, we are asking that for any linear series $\mathcal A'$ (of degree $d'$ and dimension $r'$)
contained in the linear series associated to $\psi$
the inequality ${d'}/{r'}\geq d/r$ has to be satisfied.



\begin{rem}\label{clifford}
{\upshape
It is easy to see that  when $\psi$ is induced from a line bundle, it is sufficient to check the
inequality  for any {\em complete} linear series contained in the one associated to $\psi$.
The classical results on divisors on curves, such as Clifford's Theorem and its generalizations
(\cite{Bv}, \cite{reid}) and the Riemann-Roch Theorem, imply quite easily the following results
(cf. \cite{ACG2} and \cite{tesi}).

\medskip

$(1)$  If  $C$ is a non-hyperelliptic curve, the canonical embedding of $C$ in $\pr^{g-1}$ is linearly stable.

$(2)$ If $C$ is hyperelliptic, the canonical morphism is linearly semistable, but not stable.

$(3)$ If $C$ is a smooth curve of genus $g\geq 1$ and $L$ is a line bundle on $C$ of degree $d>2g$,
the embedding induced by $L$ is linearly stable.
}
\end{rem}

We are interested in the linear stability of {\em projections} from the canonical image of a curve.

\begin{ex}
{\upshape If $C$ is a \emph{trigonal} curve, the projection from a point outside the canonical image can be linearly unstable.
Indeed, consider any effective divisor $P_1+P_2+P_3$ belonging to  the $g^1_3$ on $C$ (the $P_i$'s need not be distinct).
By the geometric  Riemann-Roch Theorem these points span a line $\ell\subset \pr^{g-1}$.
Let $P$ be a point of $\ell$ disjoint from $C$.
It can be easily checked that the projection $\pi$ from $P$ is a birational morphism.
The image of $\pi$, $\overline C$, has a triple point $R$.
If we consider the projection $\psi$ from $R$, we have, for $g\geq 5$
$$\rd(\overline C,\psi)=\frac{2g-5}{g-3}<\frac{2g-2}{g-2}=\rd(C,\pi).$$
}
\end{ex}


>From now on, $C$ is a non-hyperelliptic curve embedded in $\pr^{g-1}$ by its canonical system.

Let  $p_\Lambda$ be the projection  from a $(k-1)$-plane $\Lambda$ disjoint from $C$.
We search for conditions for  $p_\Lambda$ to be linearly stable.

Call $V=\ann (\Lambda)  \subset H^0(\om_C)$ the linear system associated to
$p_\Lambda$, and $\mathcal V=|V|$ the  associated linear series.
Let $W\subseteq V$ be any  proper subsystem. We call $\mathcal W=|W|$ the linear series,
and $\overline{\mathcal W}$ the base point free linear series obtained from $\mathcal W$ by
 eliminating the base points.

If $\dim W=g-k-\alpha$, and $\deg \overline{\mathcal W}= 2g-2-d$, then  $W$ is
{\em not}   destabilising for $V$ if and only if
$$\, \alpha \geq d\, \frac{g-k-1}{2g-2}.$$
Let $D$ be  the  effective divisor  of base points of $\mathcal W$.
Roughly speaking, the inequality above implies that $D$
should impose ``enough'' conditions on $V$ itself.
Indeed, a \emph{sufficient} condition for $W$ not to be destabilising is
\begin{equation}\label{quellavera}\dim V(-D)\leq \dim V- d \frac{g-k-1}{2g-2},\end{equation}
where, as usual, $ V(-D)=V\cap H^0(\om_C(-D))$.
The geometric meaning of  this  condition  is that the $(k-1)$-plane $ \Lambda$
intersect the $(d + h^0(D))$-plane spanned by $D$ in a plane of dimension smaller or equal to
$g- h^0(\om_C(-D))- d (g-k-1)/(2g-2) -1$.

\begin{rem}
{\upshape
Using the stronger versions of Clifford Theorem proved in \cite{Bv} and \cite{reid}, it can be shown that
the projection of a canonical non-trigonal curve from {\em any} point not contained in  it is linearly stable.
Moreover, one can show that the projection of a trigonal canonical curve from a point {\em not} contained
in a trisecant line is linearly stable (cf. \cite{tesi}).
In what follows, we generalize these results for projections from positive dimensional subspaces.
}
\end{rem}


Given a line bundle $L$ over a smooth curve $C$, its Clifford index is
$\cliff (L)=\deg L-2(h^0(L)-1)$. If $D$ is a divisor on $C$, $\cliff(D):=\cliff (\oo_C(D))$.

\begin{defi}\label{ci} The \emph{Clifford index} of a curve $C$ of genus $g\geq 4$ is the  integer:
$$\cliff (C)=\min \{\cliff (L) \,\,\,\, |\,\,\, L\in \pic(C), \,h^0(L)\geq 2 ,\, h^1(L)\geq 2 \}.$$
When $g=2$ we set $\cliff (C)=0$; when $g=3$ we set $\cliff(C)=0$ or $1$ according to whether $C$ is hyperelliptic or not.
\end{defi}
A line bundle with $h^0$ and $h^1$ greater or equal to $2$ is said to \emph{contribute} to the Clifford index.
Brill-Noether theory
shows that $\cliff (C)\leq \left[(g-1)/2\right]$, and that
equality holds if  $C$ is general in moduli. Clifford's Theorem
says that the curves with Clifford index $0$ are exactly the
hyperelliptic ones. It is easy to prove that the curves with
Clifford index $1$ are the trigonal ones and the smooth plane
quintics. In general, the Clifford index and the gonality of a
curve $C$ are related by the following (cf. \cite{C-M})
$$ \gon (C)-3\leq \cliff (C)\leq \gon (C)-2. $$

\begin{rem}\label{maledizione}\upshape{
As the Clifford index of a curve $C$ measures how large is the ratio between the degree and the
dimension of special linear series on $C$,
it seems natural to guess that the canonical curves with higher Clifford index have linearly stable
projections from positive-dimensional subspaces of $\pr^{g-1}$.
However, this guess is false.
The problem is that the Clifford index does not control the divisors having $H^1$ of dimension $1$.
Indeed, consider a non-hyperelliptic curve $C$ with \emph{arbitrary} Clifford index, and let
$D=P_0+\ldots P_k$ be an effective divisor consisting of $k+1$ points that impose independent
conditions on $H^0(\om_C)$.
Consider a section $\varphi$ of $H^0(\om_C)$ not vanishing at anyone of the $P_i$'s
(a general section will do).
The linear subsystem $V$ of $H^0(\om_C)$ spanned by $H^0(\om_C(-D))$ and by $\varphi$
has no base points by construction, and
has dimension $g-k$.
Hence, $V$ induces the projection of the canonical image of $C$ from a subspace of projective
dimension $k-1$ disjoint from it.
As soon as $k>1$ this projection is linearly unstable, because
$$
\rd (V)=\frac{2g-2}{g-k-1}>\rd (\om_C(-D))=\frac{2g-3-k}{g-k-2}.
$$
Note that $h^1(\om_C(-D))=1$, and hence $\om_C(-D)$ is one of the divisors that does not contribute
to the Clifford index of $C$.
}
\end{rem}

\begin{prop}\label{stabprojgenerale}
Let $C\subset \pr^{g-1}$ be a canonical curve, and $k$ an integer such that  $\cliff (C)\geq 2k$.
Let $\Lambda $ be a $(k-1)$-plane in $\pr ^{g-1}$ disjoint from $C$ such that
\begin{equation}\label{vabe}
\dim (\Lambda \cap \spn (D))< d\frac{g+k-1}{2g-2}-1
\end{equation}
for any special effective divisor $D$ on $C$ with degree $d\leq 2k-1$ such that $\dim \spn (D)=d-1$.

Then the projection with centre $\Lambda$ is linearly stable.
\end{prop}

\begin{proof}
Let $V\subset H^0(C,\omega_C)$ be the linear system associated to the projection with centre $\Lambda$.
Let $W\subset V$ be any linear subsystem with $\dim W\geq 2$; we need to check that
$\rd (C,W)\geq \rd(C, V)$.
Let $L_W$ be the line bundle generated by the sections of $W$.
$L_W\cong \omega_C(-D)$, with $D=\ann(W)\cap C$.
Observe  that
$$\Lambda \cap \spn (D)=\pr ( \ann (V+ H^0(\om_C(-D)))\subseteq \pr (H^0(\om_C)^\vee )=\pr ^{g-1}.$$
Applying Grassman formula to $V+H^ 0(\omega_C(-D))$, condition (\ref{quellavera}) translates as:
\begin{equation}\label{num}
\dim (\Lambda \cap \spn (D))\leq \frac{k}{2g-2}d +\frac{\cliff (D)}{2}-1.
\end{equation}
Note that $\spn(D)$ is a plane of dimension $(d+c)/2-1$.

If $D$ contributes to the Clifford index of $C$ then inequality (\ref{num}) is trivially satisfied as the
right side term is bigger than $k-1$, which is the dimension of $\Lambda$.

\smallskip

If, on the other hand, $D$ does not contribute to the Clifford index, necessarily we have $h^0(\oo_C(D))=1$, because
$h^1(\oo_C(D))=h^0(\omega_C(-D))\geq \dim W\geq 2$.
By  the geometric version of the Riemann-Roch Theorem,  the points of $D$ are in general position
(i.e. $\dim \spn (D)=d-1$).
Moreover, notice that in this case $d=\cliff (D)$.

If $d> 2k(g-1)/(g+k-1)$, then condition (\ref{num}) is satisfied with strict inequality, because the space $\spn (D)$
has dimension strictly smaller than the number on the right hand side.
Hence  we can consider the case $$d \leq\frac{2k(g-1)}{g+k-1}.$$
In particular, $d$ has to be smaller or equal to $2k-1$, and under this assumption, inequality (\ref{num})
is implied by (\ref{vabe}).
Hence, the proof is concluded.
\end{proof}


For the applications contained in the next section, we need to treat the following situation.
Suppose that we are given a linear subspace $\sig$  of $\pr^{g-1}$ of dimension $\dim \sig = s-1$ (without any assumption on it).
We want to find the biggest possible integer $k$ such that there exists a linear subspace $\ld$ of dimension $k-1$ contained in $\sig$ such that the projection $\pil$ with centre $\ld$ is linearly stable.
Of course $k$ will depend on the dimension of $\sig$.

\begin{teo}\label{propkappa}
Let $C\subset \pr^{g-1}$ be a canonical non-hyperelliptic curve. Let $\sig\subset
\pr^{g-1}$ be a proper $(s-1)$-space. Let $k$ be any
positive integer smaller or equal to $\min\left\{[s/2], [\cliff (C)/2]\right\}$.
Then there is a non-empty open set of $(k-1)$-spaces contained in
$\sig$ that induce linearly stable projections of degree $2g-2$.
\end{teo}
\begin{proof}We try and find a linear space $\Lambda\subseteq \Sigma $ satisfying
the assumptions of Proposition \ref{stabprojgenerale}.
We can  replace  conditions (\ref{vabe}) with the following (more restrictive) ones:
\begin{equation}\label{dimezzi}
\dim (\Lambda \cap \spn (D))\leq \frac{d}{2}-1
\end{equation}
for any special divisor $D$ on $C$ with degree $d\leq 2k-1$ such that $\dim \spn (D)=d-1$.

Observe that  condition (\ref{dimezzi}) for $d$ even is implied by the same condition for $d+1$.
Hence we can suppose $d$ odd.
We seek the existence of $\ld$ in $\sig$ that does not contain any $((d-1)/2)$-space contained in
the span of $d$ points in general position.
Let us bound from the above the dimension of such ``bad'' spaces in the grassmanian $\gr(k,g)$
of $(k-1)$-spaces in $\pr^{g-1}$.
\begin{itemize}
\item The dimension of the spaces $\spn (D)$ is $d$.
\item The dimension of the $((d-1)/2)$-spaces contained in a fixed $(d-1)$-space $\spn (D)$
is $\dim \gr((d+1)/2,d)=d^2/4-1/4$.
\item The dimension of the $(k-1)$-spaces contained in $\sig$ that contain a fixed
$(d-1)/2$-space is $\dim \gr(k-(d-1)/2-1,s-(d-1)/2-1)=(k-(d+1)/2)(s-k)$.
\end{itemize}
Hence there exists a $(k-1)$-space in $\sig$ satisfying conditions (\ref{dimezzi}) as soon as the
grassmanian of $(k-1)$-planes contained in $\sig$ has dimension strictly greater than the
dimension of the ``bad'' family, i.e.
$$
k(s-k)=\dim \gr (k,s)> d+\frac{d^2-1}{4}+ \left(k-\frac{d+1}{2}\right)(s-k),
$$
which becomes
\begin{equation}\label{condit}
s\geq k+1+\frac{d+1}{2}.
\end{equation}
As $d$ varies from $1$ to $2k-1$, we see that the inequality obtained is $s\geq 2k+1$.

For $d=2k-1$, we can slightly improve the bound arguing as follows.
Inequality (\ref{dimezzi}) for $d=2k-1$ means that $\ld$ in $\sig$
\emph{is not entirely contained} in any $(2k-2)$-space $\spn (D)$.
Let us make the following remark

\medskip

{\em
If $\sig$ is not the whole $\pr^{g-1}$, then for any $r\leq s$, there is at most a  finite number of $r$-secant $(r-1)$-spaces \emph{entirely} contained in $\sig$.
}

\medskip

Indeed, if there were a positive dimensional family of $d$ secant $(r-1)$-spaces contained in $\sig$,
then the whole curve $C$ would be contained in $\sig$, contradicting the fact that the canonical
morphism is non-degenerate.

Therefore, the $(k-1)$-spaces \emph{contained in $\sig$} that are also contained in a $\spn (D)$
are of dimension at most
$$
2k-1+ \dim \gr\left(k,2k-2\right)=k^2-1,
$$
and the same argument as above gives the bound
$k(s-k)> k^2-1,$
hence $s\geq 2k$.

Noting that for $d\leq 2k-3$ conditions (\ref{condit}) are satisfied for $s\geq 2k$, we can conclude the proof.
\end{proof}

\begin{rem}
\upshape{Note that the condition $\cliff(C)\geq 2k$ implies necessarily that $g$ has to be
greater or equal to $4k+1$.
Hence, if we consider for instance $\sig=\pr^{g-1}$, the above result is empty for $4k\geq g$.
However, it implies for instance that if $C$ has general Clifford index (which is a general condition)
 then there exists a linear space of dimension $\left[\frac{g-1}{4}\right]-1$ such that the projection
 from it is a linearly stable map.
It has to be remarked anyway that for $k$ ``big'' with respect to $g$, the sufficient conditions
made in the proof of Theorem \ref{propkappa} to simplify the original conditions for stability
contained in Proposition \ref{stabprojgenerale}, become consistently restrictive.}
\end{rem}

\section{Application to the slope of fibred surfaces}\label{appl}

Let $f\colon S \longrightarrow B$ be a non-Albanese fibration. We are
interested in giving a lower bound for the slope $s(f)$ as an
increasing function of $q_f$. For this we will apply relative
projections to the relative canonical map $S \dashrightarrow
\mathbb{P}(f_* \omega _f)$ which induce, on the general fibre $F$,
a linearly stable projection. The bigger the center of the
projection is, the better is the bound we get. In the analysis of
linear stability of projections of canonical curves in the
previous section, appears as a fundamental ingredient the Clifford
index. As we will see, the bound we get involves naturally this
two invariants:  the Clifford index of the general fibre and the
relative irregularity $q_f$.

Taking any linear subspace of the canonical embedding of a
concrete fibre $F$ we are not sure we can extend it to a {\it
relative} linear subspace over $B$ (in order to make a relative
projection), except it is contained in the trivial part
$\mathbb{P}({\cal O}_B^{\oplus q_f})$ of the Fujita decomposition
$$f_* \omega_f ={\cal A} \oplus {\cal O}_B^{\oplus q_f}.$$

Moreover, such an election allows us to control the degree of the
sheaves involved, since $\rm{deg} \cal A$ = ${\rm{deg}} f_*
\omega_f$.

We present here two different proofs. To the, yet classical,
method of Xiao to study the slope of fibrations, has recently
joined the generalized Cornalba-Harris method. Although they are  of different
nature, the application of both seem to give very similar results
in several situations (cf. \cite{tesi}, \cite{LS}), at least
in the case of fibred surfaces (\cite{trigonalebarjastoppino} gives
a higher dimensional example) .
It is an intriguing question
whether both methods are in fact equivalent or not.
Our aim is to
show how either method provides, in this case, exactly the same
bound, and to present this fact as an instance of this
parallelism.

As we will see, Xiao and Cornalba-Harris start from  a subsheaf of the pushforward of a line
bundle on the total space (in our, and most cases, $f_*\of$); from this, they give as an output an inequality
involving divisor classes on the base.
However, while Xiao's method needs almost no hypothesis, the one of  Cornalba-Harris requires
a GIT stability condition on the maps induced by the subsheaf on the general fibres.
Nevertheless, as hopefully the computations made here will show, the {\em linear stability}
of the maps induced on the general fibres, although not required by Xiao, is a fundamental
ingredient for both the approaches.

\bigskip

Applying our results on linear stability of projections, we are able to find a direct factor
$\mathcal E$ of $f_*\of$ which induces linearly stable projections on the general fibres,
and such that $\deg \mathcal E=\chi_f$.

\bigskip

Given a fibred surface $f\colon S\longrightarrow B$, we define its {\em  Clifford index} $\cliff (f)$ as
the maximum of the Clifford indices of the fibres (cf. \cite{konnocliff}).
As $\cliff$ is a lower semicontinuous locally constant function, $\cliff (f)$ is the Clifford index of the general fibres.

\begin{prop}\label{ole}
Let $f\colon S\longrightarrow B$ be a fibred surface.
If $k=\min \{[\cliff (f)/2], [\qb/2]\}$, there exists a decomposition
$$f_*\of=\mathcal E\oplus \oo^{\oplus k}_B$$
such that the fibre of $\mathcal E$ on general $t\in B$ is a  linear system inducing a linearly stable
degree $2g-2$ morphism
of the fibre $f^{-1}(t)=F_t$.
\end{prop}
\begin{proof}
If $f$ is an Albanese  fibration, or if $\cliff (f)\leq 1$, $\mathcal E$ is  the whole sheaf $f_*\of$,
and the statement  is satisfied,
because for a general fibre $F$, $H^0(F, \omega_F)$ is base-point-free, and it induces a linearly
stable embedding  (Remark \ref{clifford}).

Otherwise, let us consider the Fujita decomposition
$$f_*\of=\mathcal A \oplus \oo_B^{\oplus \qb}.$$

The sheaf $\mathcal A$ induces on a fibre $F_t$ a projection of
the canonical image from  the $(\qb-1)$-plane $\Sigma_t =\pr
\!(\ann (\mathcal A\otimes \kapp (t))$ (of course $\Sigma_t $ is
canonically identified with $\pr (\oo_B^{\oplus \qb}\otimes \kapp
(t))$).

A general fibre $F$  is smooth and $\cliff (F)=\cliff f$.
Let us fix such a general fibre, and drop the small $t$ from the notations.

Let $\Lambda$ be a $k-1$-plane contained in $\Sigma$  and let
$A'\subseteq H^0(\om_{F})$ be the linear system associated with
the projection from it. Note that, as $\oo_B^{\oplus \qb}$ is
trivial, we can extend $\Lambda$ to a trivial direct factor of
$f_*\of$ and get a decomposition
$$
f_*\of =\mathcal E\oplus \oo^{\oplus k}_B.
$$
By Theorem \ref{propkappa}, as conditions $2k\leq \qb$  and $\cliff(F)\geq 2k$ are
satisfied, on any general fibre there exists a dense open set of $(k-1)$-plane contained in $\Sigma$
inducing linearly stable, base-point-free  projections of degree $2g-2$.
So we can choose one $\Lambda$ in our fixed fibre $F$ such that
the fibre of the corresponding $\mathcal E$ on general  $t$ enjoys the same properties.
\end{proof}


\bigskip

We now come to the two proofs of Theorem \ref{resultado}.


\subsection*{Via Xiao's method}

Xiao's method is a well established way of studying slopes of
fibrations (cf. \cite{X}, \cite{A-K}, \cite{Barja3folds},
\cite{Ohno}, \cite{konnotrig}, \cite{konnohyp}). We just sketch it
and refer to \cite{A-K} and to \cite{X} for details. Consider the
Harder-Narashimann filtration of any subsheaf ${\cal F}$  of $ f_*
\omega_f$:

$$
0=:{\mathcal E}_0 \subset{\mathcal
E}_1\subset\cdots\subset{\mathcal E}_n={\cal F}$$

\noindent and let $\mu_1>\cdots>\mu_n$ ($\mu_i:=\mu({\cal
E}_i/{\cal E}_{i-1})$) be the associated slopes. Set $r_i={\rm
rk}{\cal E}_i$. We have
$$
\deg {\cal F}=\sum_{i=1}^n r_i(\mu_i-\mu_{i+1}),\qquad
(\mbox{where }\mu_{n+1}=0)
$$

For technical reasons, it is necessary that all the sequence of
slopes is decreasing (including $\mu_{n+1}=0$), so we need $\mu_n
\geq 0$. This is always achieved if $\cal F$ is not only a
subsheaf but also a direct summand of $f_* \omega_f$ (which is a
nef vector bundle on $B$).

For each $i$, the composite of the natural sheaf homomorphisms
$$
f^*{\mathcal E}_i\rightarrow f^*f_*\omega_f\rightarrow \omega_f
$$
\noindent induces a rational map $S \rightarrow {\mathbb
P}_B({\mathcal E}_i)$. Up to a suitable sequence of blowing-ups
$\epsilon: {\widehat S} \longrightarrow S$ (which does not modify
the general fibre $F$), the above map becomes a morphism for every
$i$. Let $M_i$ be the moving part of the pull-back of the
tautological line bundle $H_i$ on ${\mathbb P}_B({\mathcal E}_i)$. $M_i|_F$ is
a base point free linear system on $F$ which induces a
map into ${\mathbb P}^{r_i-1}$ (a fibre of ${\mathbb
P}_B({\mathcal E}_i)\rightarrow B$), of degree $d_i$.

\begin{prop} (cf. \cite{X}) \label{Xiao}
For any sequence of indices with $1\leq i_1<\cdots<i_m\leq n$ we
have

$$
M_n^2\geq \sum_{p=1}^m
(d_{i_p}+d_{i_{p+1}})(\mu_{i_p}-\mu_{i_{p+1}})
$$
\noindent
where $i_{m+1}=n+1$.
\end{prop}

\bigskip

Then, we can proceed to give a proof of
\ref{resultado}:

\begin{proof}{\it of \ref{resultado}}

Following the notations of Proposition \ref{ole}, we put ${\cal
F}={\cal E}={\cal A} \oplus {\cal O}^{\oplus (q_f -k)} $  and
apply \ref{Xiao} for the whole set of indexes $\{1,2,...,n\}$:

$$M_n^2\geq \sum_{i=1}^n
(d_{i}+d_{i+1})(\mu_{i}-\mu_{i+1}).
$$

Since by construction the linear system ${M_n}_{\vert F}$ is
linearly stable and of degree $d_n=2g-2$, then for all $i=1,..,n$
we have

$$\frac{d_i}{r_i-1} \geq \frac{d_n}{r_n-1}=\frac{2g-2}{q_f-k-1}=:\alpha.$$

Using that $r_{i+1} \geq r_i+1$ and that ${\rm deg} {\cal E}={\rm
deg } f_*\omega _f =\chi _f$ we conclude

$$M_n^2 \geq 2\alpha {\chi _f}- \alpha \mu _1.$$

On the other hand, $M_n \leq \epsilon ^* K_f$ and both are nef, so
we have
$$K_f^2 = (\epsilon ^* K_f)^2  \geq M_n^2 \geq 2\alpha {\chi _f}- \alpha \mu _1.$$

Finally, defining now ${\cal F}=f_*\omega_f$ and taking the set of
indexes $\{1,n\}$ we obtain  $$K^2_f \geq (2g-2)\mu_1$$ \noindent
which combined with the previous inequality produces the desired
result
$$K^2_f \geq 4\frac{g-1}{g-k}\chi_f.$$

\end{proof}

\

\begin{rem}\label{semistable}
{\upshape In the particular case when $\cal A$ is semistable, we can take in the previous proof ${\cal F}={\cal A}$.
Then the same argument produce

$$K^2_f \geq 2\frac{d}{g-q_f}\chi_f.$$

If, moreover, we know that ${\mathbb P}({\cal O}_B^{\oplus q_f})$  does not meet the general fibre $F$, then $d=2g-2$ and so

$$K^2_f \geq 4\frac{g-1}{g-q_f}\chi_f.$$
}

\end{rem}


\subsection*{Via a Theorem of Cornalba and Harris}

The method of Cornalba-Harris is introduced in \cite{C-H}.
Let us summarize the version for fibred surfaces
\footnote{The original Cornalba-Harris Theorem requires the assumption of {\em Hilbert}
instead of linear stability. For curves, linear stability implies Hilbert stability as proved  in \cite{ACG2} or in \cite{tesi}.
It is not known whether the converse implication holds.} following the generalization presented in  \cite{LS}.

\begin{teo}[Cornalba-Harris]\label{princ1}
Let $f\colon S\rightarrow B$ be a fibred surface.
Let $L$ be a line bundle on $S$ and $\mathcal F$ a coherent subsheaf of $f_*L$
of rank $r$ such that for general $t\in B$ the linear system
$$\mathcal F\otimes \kapp(t)\subseteq H^0(F_t,L_{|F_t})$$
induces a linear stable map. Let $\mathcal G_h$  be a coherent
subsheaf of $f_*L^{\otimes h}$ that contains the image of the
morphism
$$
\mbox{\upshape{Sym}}^h\mathcal F\longrightarrow f_*L^{\otimes h},
$$
and coincides with it at general $t$. If $N=\rk \mathcal G_h$ is
of the form $A h+ O(1)$ and $\deg \mathcal G_h$ of the form $B
h^2+ O(h)$, the following inequality holds:
\begin{equation}\label{princ1eq}
rB-A\deg \mathcal F\geq 0.
\end{equation}
\end{teo}

\medskip

Let us consider the particular case in which $\mathcal F=f_*L$ and
$\mathcal G_h=f_*L^{\otimes h}$. By the Riemann-Roch Theorem
$$
\deg f_*L^{\otimes h}= \deg f_!L^{\otimes h} + \deg
R^1f_*L^{\otimes h}=\frac{h^2}{2}L^2-\frac{h}{2}LK_f+ \deg f_*\of
+ \deg R^1f_*L^{\otimes h}.
$$
Let $d$ be the relative degree of $L$. For large enough $h$, By
Riemann-Roch on the general fibre, $N=dh-g+1$, where $g$ is the
genus of the fibration. Suppose that $\deg R^1f_*L^{\otimes h}= C
h^2+O(h)$; in this case the computation of the leading coefficient
of $\deg\mathcal G_h$ gives:
\begin{equation}\label{princ2}
rL^2+r\,C-2d\,\deg f_*L\geq 0.
\end{equation}

\begin{proof}{\it of Theorem \ref{resultado}}

Let us use Proposition \ref{ole}.
If $k=0$ ($\cliff f\leq 1$, or $f$ is an Albanese morphism), the statement of Theorem \ref{resultado} is just the slope inequality.

Otherwise, observe that the sheaf $\mathcal E$ of Proposition
\ref{ole} satisfies the assumptions of Theorem \ref{princ1}.
Consider the morphism of sheaves
$$
\sym^h \mathcal E\longrightarrow f_*\of^{\otimes h},
$$
and call $\mathcal G_h$ its image.

On general $t$, the morphism induced by $\mathcal E\otimes \kapp (t)$  has degree $2g-2$.
Moreover, we now prove  that it  is birational. Indeed, as  a consequence of Castelnuovo's bound (cf. \cite{ACG2} Exercise B-7), either the map induced by $\mathcal E\otimes \kapp (t)$ is  birational or it  factors
through a double cover over a curve of genus at most $k$.
This last case is impossible, because it would imply that
$$\cliff (f)=\cliff (F_t)\leq \gon F_t-2\leq \g\leq k,$$ contrary to the assumption. Hence,
$$
\rk\mathcal G_h=h^0(\overline F, (j^*\oo_{\pr^{g
-k-1}}(1))^{\otimes h})=(2g-2) h+O(1),
$$
where $\overline F_t$ is the image of $F_t$.
Moreover,
$$\deg \mathcal G_h\leq \deg f_*\of ^{\otimes h}$$
because $f_*\of ^{\otimes h}$ is nef (cf. \cite{Viehweg}). Hence,
the coefficient of $h^2$ in $\deg \mathcal G_h$ is smaller than
$K^2_f/2$, and inequality (\ref{princ1eq}) implies
$$
(g-k)\frac{K^2_f}{2}-(2g-2){\chi _f} \geq 0,
$$
as claimed.
\end{proof}

\begin{rem}\label{se...}
{\upshape Suppose that, under suitable assumptions, the fibre of
$\mathcal A$ itself on general $t\in B$ was a base point free
linear system of degree $d$ which induced a linear semistable
morphism. Both the Cornalba-Harris Theorem and the method of Xiao
would give as a result the inequality
$$s(f)\geq 2\frac{d}{g-\qb},$$
which coincides with the bound of Conjecture \ref{conjecture} if $d=2g-2$, that's to
say if $\mathbb{P}({\cal O}_B^{\oplus q_f})$
is disjoint from  the general fibre.
}
\end{rem}

\section{Examples}


\begin{ex} \label{doppie}
{\upshape This example is constructed in  \cite{BPhD}, sec 4.5
(see also \cite{CS}, Example 4.1). Let $\y$ and $B$ be smooth
curves. Let $\g>0$ be the genus of $\y$. We consider $B \times
\y$. Let $p_1$ and $p_2$ be the two projections, and $H_1$, $H_2$
their general fibres. Consider a smooth divisor $R\in
|2nH_1+2mH_2|$ (by Bertini's Theorem such a divisor exists, at
least for sufficiently large $n$ and $m$). Let   $\rho\colon
X\rightarrow B\times\y$ be the double cover ramified over $R$.

Consider the fibration $f:=p_1\circ\rho \colon X\rightarrow B$;
its  general fibre is a double cover of $\y$, and its genus is $g=2\g+m-1$.
A computation shows that its  slope is
$$ s(f)=4\frac{2\g+m-2}{\g+m-1}=4\frac{g-1}{g-\g}.$$
The relative irregularity is exactly $\qb= \g$. Indeed,
$$q=
h^1(X,\oo_X)=
h^1(B\times \y,\oo_{B\times \y})+h^1(B\times \y,\el^{-1})=$$
$$=b+\g +h^1(B,K_B(nP_1))+h^1(\y,K_\y(mP_2))=b+\g.$$
Hence, notice  that this  fibrations  have  slope reaching the expected bound of
Conjecture \ref{conjecture}, (regardless to the Clifford index).
Quite interestingly, these fibrations are also examples of slope minimal with
respect to the bound for double cover fibrations established in \cite{CS}.

What about the Clifford index?  In general, the gonality of the general fibre of these fibrations
is at most twice the gonality of the quotient $\y$, and so it is smaller or equal
to $\g+3$ if $\g$ is odd, and $\g+2$ if $\g$ is even.
Hence, the Clifford index of the general fibre $X_t$ is smaller or equal
to $ \g+1$ for odd $\g$ and to $\g$ for even $\g$.
Under suitable assumptions,  the Clifford index is ``almost'' $\g$, as shown by the
following standard argument.
\begin{lem}
Suppose that  $\y$ has general gonality  $\gon (\y)=\alpha =\left[
\frac{\g+3}{2}\right]$, and suppose that $\alpha$ is a prime
number. If we  choose $m\geq \g^2+\g+4$ in the above construction,
then $\cliff (X_t)\geq \g$ for odd $\g$, and $\cliff(X_t)\geq
\g-1$ for even $\g$.
\end{lem}
\begin{proof} Let  $\beta$ be the gonality of $\xt $.
We want to prove that $\beta =2\alpha$.
By definition of gonality, $\beta\leq 2\alpha$. Let us suppose then that $\beta$ is
strictly smaller than $2\alpha$.
Consider the following diagram
$$
\begin{array}{rrcll}
 & & \xt & & \\
 & {\scriptstyle\sigma_1}\, \swarrow &\, \downarrow\, {\scriptstyle\sigma} &\searrow \,\scriptstyle{\sigma_2} & \\
 \pr^1 \!\!\!&\stackrel{\beta_1}{\longleftarrow}\!\!\!\! & \overline \xt&\!\!\!\!\stackrel{\beta_2}{\longrightarrow} & \!\!\!\pr^1\\
  & {\scriptstyle\pi_1}\, \nwarrow &  \cap & \nearrow \,{\scriptstyle\pi_2} & \\
  & & \pr^1\times \pr^1 & & \\
\end{array}
$$
where  $\sigma_1$ is  a degree $\beta$ morphism, and $\sigma_2$ the
composition of the quotient morphism
$\psi\colon \xt\rightarrow \y$ with a morphism $\y\rightarrow \pr^1$ of degree $\alpha$;
the $\pi_i$ are the projections, $\sigma=\sigma_1\times\sigma_2$ and
$\overline \xt=\sigma(\xt)$.
Let $d$ be the degree of $\sigma$; clearly $d\mid 2\alpha$ and $d\mid \beta$,
hence $d=1,2,\alpha$ are the only possibilities.
By the adjunction formula,
$$
2 p_a(\overline \xt)-2  =(K_{\pr^1\times\pr^1}\overline\xt +\overline\xt ^2)=
  2\frac{(\alpha \beta)}{d^2}-\frac{\beta}{d}-\frac{2\alpha}{d}+1.
 $$
 If $d=1$, then $g\leq p_a(\overline\xt)\leq (2\alpha -1)(\beta-1)\leq (\g+2)(\g+1)$.
 Remembering that $g=2\g+m-1$, we deduce that $m$ has to be smaller or equal to $\g^2+\g+3$,
 contrary to the assumption.

 If $d=2$, then observe that by assumption $g> \g(\g+1)+1\geq 4(p_a(\overline \xt))+1$,
 then it follows from  Lemma 1.7 of \cite{CS}, that $\overline \xt$ is isomorphic to $\y$.
 But then $\beta_1$ would be a morphism from $\y$ to $\pr^1$ of degree strictly smaller than $\alpha$,
 which is a contradiction.

It remains to deal with the case $d=\alpha$.
In this case, It has to be $\overline \xt=\pr^1$, and $\sigma=\sigma_1$.
Then, one can consider the composite morphism $ \xt \stackrel{\sigma\times \psi}{\longrightarrow} \pr^1\times\y$,
which has degree either $1$ or $\alpha$ (remember that $\alpha$ is prime).
The case of degree $1$ would imply (again by adjunction)  that $g\leq (\g+3)(3\g-1)/4+1$, so it can be excluded.
The case of degree $\alpha$ would imply that $\y\cong \pr^1$, a contradiction, because we assumed that $\g=g(\y)>0$.
\end{proof}}
\end{ex}

\bigskip

\begin{ex}\label{abeliane}
{\upshape
The following construction leads in particular to examples of fibrations with $\qb=2$ and Clifford index big.
Let $S$ be an abelian surface, and let $C$ be a smooth  curve  of
genus $g\geq 3$ contained in it. By the adjunction formula, $C^2=
(CK_S+C^2)=2g-2.$ By Riemann-Roch
$$h^0(S, \oo_S(C))=\frac{(CK_S+C^2)}{2}+\chi_S=g-1\geq 2.$$
Hence we can consider an algebraic pencil (i.e. a linear series of dimension $1$) in $|C|$.
Let $\widehat S$ be the blow up of $S$ in the $2g-2$ base points.
The pencil induces   a fibration  $f\colon\widehat S \longrightarrow \pr^1$, which clearly has relative
irregularity $\qb=q(S)=2$, and whose Clifford index  is the Clifford index of $C$.

In the following we shall prove that there exist abelian surfaces containing curves of
arbitrary genus and Clifford index big.
We will use an argument suggested to us by A. Knutsen.

\medskip

Let us first recall the following definitions and results.

A line bundle $L$ on a variety $X$ is said to be {\em $k$-very ample} if for any $0$-dimensional
scheme $Z$ of length $k+1$, the restriction map
$$H^0(X,L)\longrightarrow H^0(Z, \oo_Z(L))$$
is surjective; hence, in particular, a line bundle is  $0$-very ample iff it is globally generated,
$1$-very ample iff it is very ample, $2$-very ample iff it separates tangent vectors, and so on.
If $C$ is a smooth curve then the gonality of $C$ is $k+1$ if and only if $\omega_C$ is
$(k-1)$-very ample but not $k$-ample (this is a straightforward consequence of Riemann-Roch).

If $C$ is a smooth curve contained in a smooth projective surface $S$, by adjunction
$$\omega_C \simeq (\omega_S\otimes \oo_S(C))\otimes \oo_C.$$
Hence, we derive immediately that if $\gon (C)\leq k+1$, then $ \omega_S\otimes \oo_S(C)$
is {\em not} $k$-very ample.
We will use the following result.
\begin{teo}[Bauer-Szemberg, \cite{BS}]\label{abelian}
Let $S$ be an abelian surface with Picard number $1$, and $L$ a line bundle on $S$
of type $(1,d)$, $d\geq 1$.
Then $L$ is $k$-very ample if and only if $d\geq 2k+3$.
\end{teo}
We are now ready to prove the following
\begin{lem}
Let $S$ be an abelian surface with Picard number $1$ and let $L$ be an ample line bundle of
type $(1,d)$, with $d\geq 1$.
Then if $C$ is a smooth curve of genus $g$ contained in the linear system associated to $L$,
$\gon (C)\geq \frac{d}{2}=\frac{g+1}{2}$.
\end{lem}
\begin{proof}
Let $k+1$ be the gonality of $C$. Remember that $g(C)=d-1$.
Suppose by contradiction that $k+1<d/2$. This implies that $d\geq
2k+3$, and by Theorem \ref{abelian},
$\oo_S(C)\simeq\omega_S\otimes \oo_S(C)$ is $k$-very ample. From
the above remarks it follows that $\gon (C)>k+1$, which is the
desired contradiction.
\end{proof}
It is worth noticing that the construction above can be made in much more generality using the results of \cite{terakawa}.

Hence, we can construct fibrations from an Abelian surface to $\pr^1$ with ``almost general'' Clifford index.
\begin{rem}\label{controesempio}
\upshape{
Note that these fibrations all  have slope $6$.
Indeed, given any such fibration $f\colon \widehat S\longrightarrow \pr^1$,
$$K_f= K_{\widehat S}-f^*K_{\pr^1}\sim \sum_{i=1}^{2g-2}E_i+2C,$$
where $E_i$ are the exceptional divisors of the blow up $\widehat S\rightarrow S$.
Hence,
$$K^{2}_f=\left(\sum_{i=1}^{2g-2}E_i+2C\cdot \sum_{i=1}^{2g-2}E_i+2C\right)=
4\left(\sum_{i=1}^{2g-2}E_1\cdot C\right)+\left(\sum_{i=1}^{2g-2}E_i\cdot \sum_{i=1}^{2g-2}E_i\right)=6(g-1),$$
and
$$\deg f_*\of =\chi_{\widehat S}-\chi_{\pr^1}\chi_F= g-1.$$
This slope is coherent, and indeed bigger than, the bound given by Theorem \ref{resultado}, which is $4$.

It  is also coherent with the bound of Conjecture
\ref{conjecture}, for any genus except for
$g=3$, when it gives a counterexample for the case $\qb=g-1$.
}
\end{rem}

\begin{rem}{\upshape
One could make an analogous construction starting from a $K3$
surface. By a result of Knutsen (\cite{K3K}) there are $K3$
surfaces containing curves of any possible gonality. Hence this
construction leads to fibrations with $\qb=0$ and $\cliff (f)$
arbitrary. In this case the slope is $6\frac{g-1}{g+1}$. Note that
this slope reaches  exactly the bound for fibrations with general
Clifford index  and odd genus found by Konno (cf.
\cite{konnocliff}, \cite{A-K}) and by Eisenbud-Harris for
semistable fibrations(\cite{Ha-Ei}).}
\end{rem}
}
\end{ex}


\bigskip
\bigskip
\noindent Miguel \'Angel Barja,  Departament de Matem\`atica  Aplicada I, Universitat Polit\`ecnica de Catalunya, ETSEIB Avda. Diagonal, 08028 Barcelona (Spain).
E-mail: \textsl{Miguel.Angel.Barja@upc.edu}

\bigskip
\noindent Lidia Stoppino, Dipartimento di Matematica, Universit\`a di Roma TRE, Largo S. L. Murialdo, 1  I-00146, Roma  (Italy).
E-mail: \textsl {stoppino@mat.uniroma3.it}.

\end{document}